\newcommand{\x}{\mathbf{x}}
\newcommand{\y}{\mathbf{y}}
\newcommand{\w}{\mathbf{w}}
\newcommand{\z}{\mathbf{z}}
\newcommand{\e}{\mathbf{e}}
\newcommand{\F}{\mathbf{F}}
\newcommand{\bP}{\mathbf{P}}
\newcommand{\bE}{\mathbf{E}}

\newtheorem{Assumption}{Assumption}

\documentclass[twocolumn]{autart}    

\usepackage{amsmath,amssymb,amsfonts}
\usepackage{graphicx}          
\usepackage{xcolor}
\usepackage{enumitem}
\usepackage{multirow}
\usepackage{float}
\begin{document}

\begin{frontmatter}

\title{Explicit construction of 
the minimum error variance estimator for stochastic LTI-ss systems.\thanksref{footnoteinfo}} 

\thanks[footnoteinfo]{Corresponding author D.~Eringis.}

\author[aauES]{Deividas Eringis}\ead{der@es.aau.dk},    
\author[aauES]{John Leth}\ead{jjl@es.aau.dk},               
\author[aauES]{Zheng-Hua Tan}\ead{zt@es.aau.dk},  
\author[aauES]{Rafal Wisniewski}\ead{raf@es.aau.dk}, 
\author[Fr]{Mihaly Petreczky}\ead{mihaly.petreczky@centralelille.fr}

\address[aauES]{Dept. of Electronic Systems, Aalborg University, Aalborg, Denmark}  
\address[Fr]{Laboratoire Signal et Automatique de Lille (CRIStAL), Lille, France}             


\begin{abstract}                          
    We showcase the derivation of the optimal (minimum error variance) estimator, when one part of the stochastic LTI system outputs is not measured but is able to be predicted from the measured system outputs. 
\end{abstract}

\end{frontmatter}

\section{Introduction}
Realization theory of stochastic linear time invariant state-space
representations (LTI-ss) with exogenous inputs is a mature
theory \cite{LindquistBook,Katayama:05}. In particular,
there is constructive theory for a 
minimal stochastic LTI-ss representation of a process $\y$ with exogenous input $\w$. The construction uses geometric ideas, and it is based on oblique projection of future outputs onto past inputs and outputs.\\
    Note that, in system identification it is
often assumed that jointly $(\y,\w)$ has a realization by an autonomous stochastic LTI system driven by white noise.
Indeed, if $\w$ has a realization by a stochastic LTI-ss representation driven by i.i.d gaussian noise, and $\y$ has a realization by a LTI-ss representation with exogenous input
$\w$ and i.i.d. gaussian noise, then under some mild assumptions 
(absence of feedback from $\y$ to $\w$)
$(\y,\w)$ will be the output of an autonomous stochastic LTI-ss representation.  It is then natural 
to ask the question how to construct a minimal stochastic
LTI-ss realization of $\y$ with input $\w$, from an LTI-ss realization of the joint process $(\y,\w)$, instead of
computing a realization of $\y$ using oblique projections. \\
    In this paper we present an explicit construction of a minimal stochastic LTI-ss representation
of $\y$ with an exogenous input $\w$ from an autonomous stochastic LTI-ss representation of the joint process $(\y,\w)$.  The basic idea is as follows: we will assume that $(\y,\w)$ is stationary, square-integrable, zero-mean, jointly Gaussian stochastic processes and there is no feedback from $\y$ to $\w$. Then use
the result of \cite{jozsa2018relationship} stating that
there exists a minimal LTI-ss realization of $(\y,\w)$
with matrices which admit a upper-triangular form. This allowed us to separate out part of the innovation noise of $(\y,\w)$, which purely drives $\w$, thus allowing us to formulate this construction.\\
    Our motivation for developing an explicit construction of
an LTI-ss realization of $\y$ with input $\w$ from a LTI-ss realization of $(\y,\w)$ was that 
this construction turned out to be useful in deriving non-asymptotic error bounds of PAC-Bayesian type \cite{alquier-15} for LTI-ss systems \cite{CDC21paper}. The latter could be a first step towards extending the PAC-Bayesian framework for stochastic state-space representations. \\
    More precisely, one of the byproducts of the construction of this paper is a one-to-one relationship between LTI-ss
systems which generate $(\y,\w)$ and optimal linear estimators of future values of $\y$ based on past values of $\w$. This relationship is useful in Bayesian learning algorithms, when one needs to define a parameterised set of predictors (Hypothesis class). All prior knowledge or uncertainty in the data generating system can then easily be mapped to knowledge or uncertainty of the predictor. \\
    The contribution of the paper can also be viewed as 
as follows.
We wish to construct an estimator of $\y(t)$ given past $(s<t)$ and present $(s=t)$ measurements of $\w(s)$. We consider a specific class of relationships, specifically when the two processes are related by a common stochastic LTI-ss system, i.e., $[\y^T(t)\;\w^T(t)]^T$ is an output of an LTI-ss system. The problem of finding this estimator can also be thought of as trying to estimate non-measurable quantities of a system from measurable quantities. 

\textbf{Related work:}
As it was pointed out above, stochastic realization theory
with inputs is a mature topic with several publications,
see the monographs \cite{LindquistBook,Katayama:05,CainesBook} and
the references therein. However, we have not found in the literature an explicit procedure for constructing a  stochastic LTI-ss realization in forward innovation form of $\y$ with input $\w$ from the joint stochastic LTI-ss realization of $(\y,\w)$. The current note is intended to fill this gap. \\
    We will need to further analyse the relationship between $\y$ and $\w$ by feedback-free assumption. In \cite{GRANGER196328}, the author defines what it means for one process to cause another, a similar notion to feedback. In \cite{CainesFeedback}, the authors further extend the notion and define weak and strong feedback free processes. As strong feedback free condition implies weak feedback free, we consider the relaxed case of weak feedback free throughout the paper. In frequency domain using causal real rational transfer function matrices to describe processes $\y$ and $\w$, and analysing these processes with feedback free assumption, yields a straightforward construction of estimator of $\y$ given $\w$, see \cite{feedbackProcesses} and \cite{Gevers1982OnJS}.
In this paper we study this problem in time domain, using LTI-ss representations. 

\textbf{Outline:}
This paper is organised as follows. Below we start by defining the notation and terminology used in this paper, then in Section \ref{sec:sec1} we reformulate the state-space system driven by innovation of $[\y^T(t)\;\w^T(t)]^T$ into a state-space system, which yields a realisation of $\y$, driven by $\w$ and the innovation of a purely non-deterministic part of $\y$. Afterwards in Section \ref{sec:OptimalPred}, given this new realisation we provide the optimal (in the sense of minimum error variance) estimate of $\y$.

{\bf Notation and terminology} Let $\F$ denote a $\sigma$-algebra on the set $\Omega$ and $\bP$ be a probability measure on $\F$. Unless otherwise stated all probabilistic considerations will be with respect to the probability space $(\Omega,\F,\bP)$. In this paragraph let  $\mathbb{E}$ denote some euclidean space. We associate with $\mathbb{E}$ the topology generated by the 2-norm $||\cdot||_2$, and the Borel $\sigma$-algebra generated by the open sets of $\mathbb{E}$. The closure of a set $M$ is denoted $clM$. For $S\subseteq\mathbb{N}$ and stochastic variables $\y,\z_1,\z_2,\dots$ with values in $\mathbb{R}$ we denote by $\bE(\y~|~\{\z_i\}_{i\in S})$ the conditional expectation of $\y$ with respect to the $\sigma$-algebra $\sigma(\{\z_i\})$ generated by the family $\{\z_i\}_{i\in S}$.  Recall that $\bE(\z\x)$ define an inner product in $L^2(\Omega,\F,\bP)$ and that $\bE(\y~|~\{\z_i\}_{i\in S})$ can be interpreted as the orthogonal projection onto the closed subspace $L^2(\Omega,\sigma(\{\z_i\}_{i\in S}),\bP)$ which also can be identified with the closure of the subspace generated by $\{\z_i\}_{i\in S}$. That is, 
\begin{align}\label{zigen}
	\textstyle
	L^2(\Omega,\sigma(\{\z_i\}_{i\in S}),\bP)=cl\left\{\sum_{i\in S}\alpha_i\z_i~|~\alpha_i\in\mathbb{R}\right\}
\end{align}
with only a finite number of summands in \eqref{zigen} being nonzero when $S=\mathbb{N}$. Moreover, for a closed subspace $H$ of $L^2(\Omega,\F,\bP)$ and a stochastic variable $\y$ with values in $\mathbb{E}$ and $\bE(||\y||_2^2)<\infty$, we let $\bE(\y~|~H)$ denote the $\dim(\mathbb{E})$-dimensional vector with $i$th coordinate equal to $\bE(\y_i~|~H)$ with $\y_i$ denoting the $i$th coordinate of $\y$. 

There are two closed subspaces of particular importance. Following \cite{LindquistBook}, for a discrete time stochastic process $\z(t)$ with values in $\mathbb{E}$ and $\bE(||z(t)||_2^2)<\infty$, we write $H_t^-(\z)$ for the closure of the subspace in $L^2(\Omega,\F,\bP)$ generated by the coordinate functions $\z_i(s)$ of $\z(s)$ for all $s<t$. That is, 
\begin{align}\label{h-}
	\textstyle
	H^-_t(\z)=cl\left\{\sum_{i=-\infty}^{t-1}\alpha_i^T\z(i)~|~\alpha_i\in\mathbb{E}\right\}
\end{align}
with ${}^T$ indicating transpose and only a finite number of summands in \eqref{h-} being nonzero. In a similar manner we define 
\begin{align}
	H^+_t(\z)=&\textstyle
	cl\left\{\sum_{i=t}^{\infty}\alpha_i^T\z(i)~|~\alpha_i\in\mathbb{E}\right\},\\
	H(\z)=&\textstyle
	cl\left\{\sum_{i=-\infty}^{\infty}\alpha_i^T\z(i)~|~\alpha_i\in\mathbb{E}\right\}.
\end{align}   

Let $A$, $B$ and $C$ be closed subspaces of $L^2(\Omega,\F,\bP)$. We then define
\begin{align}
	A\vee B=cl\{a+b~|~a\in A,~b\in B\}
\end{align}
and say that $A$ and $B$ are orthogonal given $C$, denoted $A\perp B~|~C$, if
\begin{align}
	\bE\Big(\big(a-\bE(a~|~C)\big)\big(b-\bE(b~|~C)\big)\Big)=0
\end{align}
for all $a\in A$ and $b\in B$. 

We use the following notation, $\mathcal{Y}=\mathbb{R}^{p}$, $\mathcal{W}=\mathbb{R}^{q}$ and for the disjoint union $\mathcal{W}^{*}=\bigsqcup_{k=1}^{\infty} \mathcal{W}^k$ we write $w=(w_1,\ldots,w_k)$ in place of the more correct $(w,k)=((w_1,\ldots,w_k),k)$ for an element in $\mathcal{W}^{*}$. 


\section{Assumptions}\label{sec:sec1}
Suppose we want to construct an estimator of the output stochastic process $\y(t):\Omega\to\mathcal{Y}$ given a sequence of measurements as inputs obtained from the stochastic process $\w(t):\Omega\to\mathcal{W}$.
In order to narrow down and formally describe the estimation problem, we assume that the processes $\y(t)$ and $\w(t)$ can be represented as outputs of an LTI system in forward innovation form: 
\begin{Assumption}\label{as:generator}
	The processes $\y(t)$ and $\w(t)$ can be generated by a stochastic discrete-time minimal LTI system on the form
	\begin{subequations}\label{eq:generator}
		\begin{align}
			\x(t+1) & =A_g\x(t)+K_g\mathbf{e}_g(t)\\
			\begin{bmatrix}
				\y(t)\\
				\w(t)
			\end{bmatrix}
			& =C_g\x(t)+\mathbf{e}_g(t), \quad Q=\bE[\mathbf{e}_g^T(t)\mathbf{e}_g(t)]
		\end{align}
	\end{subequations}
	where $A_g \in \mathbb{R}^{n \times n},K_g \in \mathbb{R}^{n \times m},C_g=[C_\y^T,C_\w^T]^T \in \mathbb{R}^{(p+q) \times n}$ for $n \ge 0$, $m,p>0$ and $\x\in\mathbb{R}^n$, $\y\in\mathbb{R}^p$,$\w\in\mathbb{R}^q$ and $\mathbf{e}_g$ are stationary, square-integrable, zero-mean, and jointly Gaussian stochastic processes. The processes $\x$ and $\mathbf{e}_g$ are called state and noise process, respectively. Recall, that stationarity and square-integrability imply constant expectation and that the covariance matrix\\ $Cov(\y(t),\y(s)) = \bE [(\y(t)-\bE[\y(t)])(\y(s)-\bE[\y(s)])^T]$\\ only depends on time lag $(t-s)$. 
	Furthermore, we require that $A_g$ is stable (all its eigenvalues are inside the open unit circle) and that for any $t,k \in \mathbb{Z}$, $k \geq 0$, $E[\mathbf{e}_g(t)\mathbf{e}_g^T(t\!-\! k\!-\! 1)]=0$, $E[\mathbf{e}_g(t)\x^T(t-k)]=0$, i.e., the stationary Gaussian process $\mathbf{e}_g(t)$ is white noise and uncorrelated with $\x(t-k)$. 
	We identify the system \eqref{eq:generator} with the tuple
	$(A_g,K_g,C_g,I,\mathbf{e}_g)$; note that the state process $\x$ is uniquely defined by
	the infinite sum $\x(t)=\sum_{k=1}^{\infty} A_g^{k-1}K_g\mathbf{e}_g(t-k)$. 
\end{Assumption}
Before we can continue we have to consider the relationship between $\y$ and $\w$. For technical reasons we can not have feedback from $\y$ to $\w$, as $\w$ would then be determined by a dynamical relation involving the past of the process $\y$. As such we have Assumption \ref{as:NoFeedback}
\begin{Assumption}\label{as:NoFeedback}
	There is no feedback from $\y$ to $\w$, following definition 17.1.1. from \cite{LindquistBook}, i.e., $$H_t^-(\y) \perp H_t^+(\w) \mid H_t^-(\w)$$ holds, i.e., the future of $\w$ is conditionally uncorrelated with the past of $\y$, given the past of $\w$.
\end{Assumption}
As a passing remark, the no feedback assumption is equivalent to weak feedback free assumption \cite{CainesFeedback} or Granger non-causality \cite{GRANGER196328}. Thus the no feedback assumption can be stated as $\y$ does not Granger cause $\w$.
\section{Result}\label{sec:Result}
Under assumption \ref{as:NoFeedback}, there exists a similarity transformation $T$ of \eqref{eq:generator} such that $\bar{A}_g=TA_gT^{-1}$, $\bar{K}_g=TK_g$ and $\bar{C}_g=C_gT^{-1}$ are upper block triangular, specifically \eqref{eq:generator} can be represented as
\begin{subequations}\label{eq:similar}
\begin{align} 
    \begin{bmatrix}
    	\bar{\x}_{1}(t+1)\\
    	\bar{\x}_{2}(t+1)
    \end{bmatrix}\hspace{-2pt}&=\hspace{-2pt}\begin{bmatrix}A_{1,1}&A_{1,2}\\0&A_{2,2}  \end{bmatrix}\hspace{-2pt}\begin{bmatrix}
    	\bar{\x}_{1}(t)\\
    	\bar{\x}_{2}(t)
    \end{bmatrix}\hspace{-2pt}  +\hspace{-2pt}\begin{bmatrix} K_{1,1} & K_{1,2} \\ 0 & K_{2,2} \end{bmatrix}\hspace{-2pt}\begin{bmatrix} \mathbf{e}_1(t)\\ \mathbf{e}_2(t)\end{bmatrix}\\
    \begin{bmatrix} \y(t)\\ \w(t)\end{bmatrix} &= \begin{bmatrix} C_{1,1} & C_{1,2} \\ 0 & C_{2,2} \end{bmatrix}\begin{bmatrix}
    	\bar{\x}_{1}(t)\\
    	\bar{\x}_{2}(t)
    \end{bmatrix} + \begin{bmatrix} \mathbf{e}_1(t)\\ \mathbf{e}_2(t)\end{bmatrix} \label{eq:similarOuput}
\end{align}
\end{subequations}
where $[\mathbf{e}_1^T(t)\;\mathbf{e}_2^T(t)]^T=\mathbf{e}_g(t)$, and such that $(A_{2,2},C_{2,2})$ is observable. Moreover, $A_{i,j}\in\mathbb{R}^{p_i\times p_j}$, $K_{i,j}\in\mathbb{R}^{p_i\times r_j}$,$C_{i,j}\in\mathbb{R}^{r_i\times p_j}$, with $r_1=p$ and $r_2=q$. \par
\textbf{The optimal estimate $\hat{\y}(t)=\bE[\y(t)\mid H_{t+1}^-(\w)]$, in the least square sense,} is then given as the output of the following LTI system  
\begin{subequations}\label{eq:T0sys}
	\begin{align}
		\hat{x}(t+1)&=\tilde{A}\hat{x}(t)+\tilde{K}\w(t) \\
	\hat{\y}(t) &= \tilde{C}\hat{x}(t)-D_0\w(t),\\
    \tilde{A}=&\begin{bmatrix}A_{1,1}&A_{1,2}-(K_{1,2}+K_{1,1}D_0)C_{2,2}\\0&A_{2,2}-K_{2,2}C_{2,2}  \end{bmatrix},\\
    \tilde{K}=&\begin{bmatrix}K_{1,2}+K_{1,1}D_0\\K_{2,2}\end{bmatrix},\\
    \tilde{C}=&\begin{bmatrix} C_{1,1} & C_{1,2}-D_0C_{2,2} \end{bmatrix},\quad 
    D_0=Q_{1,2}Q_{2,2}^{-1}.
\end{align}
\end{subequations}
with the covariance $Q=\bE[\e_g^T(t)\e_g(t)]$ partitioned according to \eqref{eq:similar}.
\section{Derivation}
Several results can be deduced from Assumption \ref{as:NoFeedback}. First, by \cite[Proposition 2.4.2]{LindquistBook}, we obtain the following relation between projections   
\begin{align}
	E[\y(t)|H(\w)]&=E[\y(t)|H_{t+1}^-(\w)],\\ 
	E[\w(t)|H_{t}^-(\w)\vee H_{t}^-(\y)]&=E[\w(t)|H_{t}^-(\w)].  
\end{align}

Secondly, from \cite[Ch. 17]{LindquistBook} it follows that the process $\y$ can then be decomposed into a deterministic part $\y_d$ and a stochastic part $\y_s$, as follows 
\begin{align}
	\y(t)&=\y_d(t)+\y_s(t), \label{eq:decomp}\\
	\y_d(t)&=E[\y(t) | H(\w)]=E[\y(t) | H_{t+1}^-(\w)],\label{eq:ydDef}\\
	\y_s(t)&=\y(t)-\y_d(t).\label{eq:ysDef}
\end{align}
Note that, as a consequence of \eqref{eq:ydDef} and \eqref{eq:ysDef} 
$$E[\y_d(t)\y_s^T(\tau)]=0\;\forall t,\tau\; ,$$
i.e., $\y_d$ and $\y_s$ are uncorrelated.
Moreover, the process $\y_s$ can be realised by a state-space system in forward innovation form
\begin{subequations}
	\begin{align}
		\x_s(t+1)&=A_s\x_s(t)+K_s\mathbf{e}_s(t),\\
		\y_s(t)&=C_s\x_s(t)+\mathbf{e}_s(t),\\
		\mathbf{e}_s&=\y_s(t)-E[\y_s(t)|H_t^-(\y_s)].\label{esdef}
	\end{align}
\end{subequations}
Finally, from \cite[Proposition 17.1.3.]{LindquistBook} we get
\begin{equation}
	H_{t}^-(\y)\vee H_{t+1}^-(\w)=H_{t}^-(\y_s)\oplus H_{t+1}^-(\w),
\end{equation}
where $\oplus$ denotes orthogonal sum, and
\begin{equation}
	\mathbf{e}_s(t)=\y(t)-E[\y(t)|H_{t}^-(\y)\vee H_{t+1}^-(\w)]. \label{eq:es}
\end{equation}


Now consider a similarity transformation $T$ of \eqref{eq:generator} such that $\bar{A}_g=TA_gT^{-1}$, $\bar{K}_g=TK_g$ and $\bar{C}_g=C_gT^{-1}$ are upper block triangular, see \eqref{eq:similar}.
From \cite{jozsa2018relationship} it then follows that $(A_{2,2},K_{2,2},C_{2,2},\mathbf{e}_2)$ is a minimal Kalman representation of $\w$ hence $\mathbf{e}_2(t)$ is the innovation process of $\w$ i.e.,
\begin{align}
	\mathbf{e}_2(t)&=\w(t) - E[\w(t) \mid H_{t}^-(\w)]\nonumber\\
	&=\w(t) - E[\w(t) \mid H_{t}^-(\y) \vee H_{t}^-(\w)]. \label{def:e2}
\end{align}
Moreover, the transformed system \eqref{eq:similar} induces a relation between the output $\y$ and input $\w$. In detail, from \eqref{eq:similarOuput} we also have
\begin{equation}
	\mathbf{e}_2(t)=\w(t)-C_{2,2}\bar{\x}_{2}(t). \label{eq:e2}
\end{equation}
Hence, substituting \eqref{eq:e2} in \eqref{eq:similar} yields the following realisation of $\y$ 
{\normalsize
\begin{subequations}\label{inout}
	\begin{flalign}
		\begin{bmatrix}
			\bar{\x}_{1}(t+1)\\
			\bar{\x}_{2}(t+1)
		\end{bmatrix}&=\begin{bmatrix}A_{1,1}&A_{1,2}-K_{1,2}C_{2,2}\\0&A_{2,2}-K_{2,2}C_{2,2}  \end{bmatrix}\begin{bmatrix}
			\bar{\x}_{1}(t)\\
			\bar{\x}_{2}(t)
		\end{bmatrix} &&\\
		&\quad +\begin{bmatrix}K_{1,2}\\K_{2,2}\end{bmatrix}\w(t) +\begin{bmatrix} K_{1,1} \\ 0 \end{bmatrix}\mathbf{e}_1(t) &&\\
		\y(t) &= \begin{bmatrix} C_{1,1} & C_{1,2} \end{bmatrix}\begin{bmatrix}
			\bar{\x}_{1}(t)\\
			\bar{\x}_{2}(t)
		\end{bmatrix} + \mathbf{e}_1(t)&&
	\end{flalign}
\end{subequations} }
Note that $\mathbf{e}_1(t)$ is the innovation process of $\y$ (with respect to $\w$), i.e.,
\begin{equation}\label{e1}
	\mathbf{e}_1(t)=\y(t)-E[\y(t) \mid H_t^-(\y) \vee H_t^-(\w)]. 
\end{equation}
\subsection{Optimal estimate} \label{sec:OptimalPred}
The goal in this section is to derive an optimal estimate (in the sense of minimum error variance). Firstly, we claim that
\begin{equation}
	\mathbf{e}_s(t)=\mathbf{e}_1(t)-E[\y(t)|\mathbf{e}_2(t)]=\mathbf{e}_1(t)-D_0\mathbf{e}_2(t) \label{eq:esClaim}
\end{equation}
where\footnote{In order to numerically compute $D_0$, we can use \eqref{e1} to replace $\y$ with $\mathbf{e}_1(t)+\bE[\y(t) \mid H_t^-(\y) \vee H_t^-(\w)]$, and since $\mathbf{e}_2 \perp H_{t}^-(\w) \vee H_{t}^-(\y)$, we get $\bE[\bE[\y(t) \mid H_t^-(\y) \vee H_t^-(\w)]\mathbf{e}_2^T]=0$. Therefore $\bE[\y(t)\mathbf{e}_2^T(t)]=\bE[\mathbf{e}_1(t)\mathbf{e}_2^T(t)]$. In summary, one can compute $D_0$ directly from the covariance of innovation noise, i.e., $D_0=Q_{1,2}Q_{2,2}^{-1}$.} $D_0=(E[\y(t)\mathbf{e}_2^T(t)])^T(E[\mathbf{e}_2(t)\mathbf{e}_2^T(t)])^{-1}$ is the minimum variance linear estimator of $\y(t)$ given $\mathbf{e}_2(t)$, see \cite[Proposition 2.2.3.]{LindquistBook}. 
In order to show \eqref{eq:esClaim}, we first demonstrate that
\begin{equation}
	H_t^-(\y) \vee H_{t+1}^-(\w) = ( H_t^-(\y) \vee H_t^-(\w) ) \oplus H(\mathbf{e}_2(t) ), \label{orthSum}
\end{equation}
where $H(\mathbf{e}_2(t))=\{\alpha^T \mathbf{e}_2(t) \mid \alpha \in \mathbb{R}^q\}$, is the space spanned by innovation process $\mathbf{e}_2(t)$, considered only at the time $t$. By definition we have
\begin{align}
	&( H_t^-(\y) \vee H_t^-(\w) ) \vee H(\mathbf{e}_2(t) ) =\nonumber\\ 
	&\quad\qquad cl\Big\{ \sum_{i=-\infty}^{t-1} \gamma_i^T \y(i) + \sum_{i=-\infty}^{t-1} \eta_i^T \w(i)\nonumber\\ 
	&\qquad\qquad + \lambda_t^T\mathbf{e}_2(t)
 \mid \gamma_i \in \mathbb{R}^p, \eta_i\in \mathbb{R}^q, \lambda_t \in \mathbb{R}^q\Big\} \label{eq:esfe12}
\end{align}
However, using definition of $\mathbf{e}_2(t)$ from \eqref{def:e2} we have
\begin{align}
	&( H_t^-(\y) \vee H_t^-(\w) ) \vee H(\mathbf{e}_2(t) ) =\nonumber \\
	&\quad\qquad cl\Big\{ \sum_{i=-\infty}^{t-1} \gamma_i^T \y(i) + \sum_{i=-\infty}^{t-1} \eta_i^T \w(i)\nonumber\\
	&\qquad\qquad + \lambda_t^T\w(t) 
	\mid \gamma_i \in \mathbb{R}^p, \eta_i\in \mathbb{R}^q, \lambda_t \in \mathbb{R}^q \Big\},
\end{align}
which equals $H_t^-(\y)\vee H_{t+1}^-(\w)$ and therefore
\begin{align}
	( H_t^-(\y) \vee H_t^-(\w) ) \vee H(\mathbf{e}_2(t) ) =  H_t^-(\y) \vee H_{t+1}^-(\w)
\end{align}
Again from \eqref{def:e2} it follows that $\mathbf{e}_2(t) \perp H_{t}^-(\w) \vee H_{t}^-(\y)$, thus \eqref{orthSum} holds. The relation \eqref{eq:esClaim} now follows since 
\begin{align}
	&E[\y(t) \mid H_t^-(\y)\vee H_{t+1}^-(\w)]\nonumber\\
	&\qquad=E[\y(t) \mid H_t^-(\y)\vee H_t^-(\w)]
	+E[\y(t) \mid \mathbf{e}_2(t)] && \nonumber \\
	&\qquad=E[\y(t) \mid H_t^-(\y)\vee H_t^-(\w)]+D_0\mathbf{e}_2(t),
\end{align} 
and therefore, using \eqref{e1} we can see that
\begin{align}
	\mathbf{e}_s(t)
	&=\y(t)-E[\y(t)\mid H_t^-(\y) \vee H_t^-(\w)]-D_0\mathbf{e}_2(t)\\
	&=\mathbf{e}_1(t)-D_0\mathbf{e}_2(t)
\end{align}

Now from \eqref{eq:esClaim} and \eqref{eq:e2} we get
\begin{equation}
	\mathbf{e}_1(t)=\mathbf{e}_s(t)+D_0\w(t)-D_0C_{2,2}\bar{\x}_{2}(t),
\end{equation}
which can be applied to \eqref{inout} to obtain the following realization of $\y$
\begin{subequations} \label{eq:y_real_with_w}
	\begin{align}
		\begin{bmatrix}
			\bar{\x}_{1}(t+1)\\
			\bar{\x}_{2}(t+1)
		\end{bmatrix}&=\tilde{A}\begin{bmatrix}
			\bar{\x}_{1}(t)\\
			\bar{\x}_{2}(t)
		\end{bmatrix}+\tilde{K}\w(t)+\begin{bmatrix} K_{1,1} \\ 0 \end{bmatrix}\mathbf{e}_s(t),\label{xbar} \\
		\y(t) &= \tilde{C}\begin{bmatrix}
			\bar{\x}_{1}(t)\\
			\bar{\x}_{2}(t)
		\end{bmatrix} +D_0\w(t) + \mathbf{e}_s(t), \label{y}
	\end{align}
\end{subequations}
with $(\tilde{A},\tilde{K},\tilde{C},D_0)$ according to \eqref{eq:T0sys}.
Finally we are in a position to derive a formula for the minimum error variance estimate $E[\y(t)\mid H_{t+1}^-(\w)]$. That is, a formula for the orthogonal projection of $\y(t)$ given past and present values of $\w$. First define $\hat{x}_g(t)=E[\bar{\x}(t)\mid H_{t+1}^-(\w)]$, then from \eqref{y} we get 
\begin{align}
	&E[\y(t)\mid H_{t+1}^-(\w)]\nonumber\\
	&\hspace{10mm}= E[\tilde{C}\bar{\x}(t)+D_0\w(t) + \mathbf{e}_s(t)\mid H_{t+1}^-(\w)]\\
	&\hspace{10mm}=\tilde{C} \hat{x}_g(t) +D_0\w(t) + E[\mathbf{e}_s(t)|H_{t+1}^-(\w)]\\
	&\hspace{10mm}=\tilde{C} \hat{x}_g(t) +D_0\w(t)\label{last}
\end{align}
where \eqref{last} follows from \eqref{eq:es}.
Now \eqref{xbar} can be used to derive a dynamical expression for $\hat{x}_g$ as follows
\begin{flalign}
	&E[\bar{\x}(t+1)\mid H_{t+2}^-(\w)] && \nonumber \\
	&=E\Bigg [ \tilde{A}\bar{\x}(t)+\tilde{K}\w(t)+\begin{bmatrix} K_{1,1} \\ 0 \end{bmatrix}\mathbf{e}_s(t) \Bigg | H_{t+2}^-(\w) \Bigg ] && \label{proxbar}
\end{flalign}
Clearly $E[\w(t)|H_{t+2}^-(\w)]=\w(t)$. For the state projection in \eqref{proxbar} we have $E[\bar{\x}(t)|H_{t+2}^-(\w)]=E[\bar{\x}(t)|H_{t+1}^-(\w)]$ since the state vector $\bar{\x}(t)$ can be expressed as an infinite sum using \eqref{xbar}, where \eqref{eq:es} is used for $\mathbf{e}_s(t)$
\begin{multline}
	\bar{\x}(t)
	=\sum_{i=1}^\infty M_i\w(t-i) + \sum_{i=1}^\infty N_i\y(t-i)\\ + \sum_{i=1}^\infty N_iE[\y(t-i)\mid H_{t-i}^-(\y) \vee H_{t-i+1}^-(\w)]\label{eq:barxsum}
\end{multline}
Finally, from \eqref{esdef} we observe that  
\begin{align*}
	E[\mathbf{e}_s(t)|H_{t+2}^-(\w)] &= E[\y_s(t)|H_{t+2}^-(\w)]\\
	&-E[E[\y_s(t)|H_t^-(\y_s)]|H_{t+2}^-(\w)]=0 
\end{align*}
since $H(\y_s)\perp H(\w)$ by \eqref{eq:ysDef}.
Finally we have obtained \eqref{eq:T0sys} the formula for the minimum prediction error variance estimate of $\y(t)$ based on $H_{t+1}^-(\w)$ (present and past of process $\w$).


An explicit construction of a realization of $\y$ with
input $\w$ has the potential to provide
an alternative to existing system identification algorithms. There are many subtleties in the consistency
analysis of subspace identification algorithms with inputs
\cite{LindquistBook,Katayama:05,chiuso2004ill}, so an alternative approach involving the identification of
an autonomous model of $(\y,\w)$ could be advantageous in some cases. \\

Before moving on to the numerical example we mention that the proposed construction is very useful when trying to learn estimators from data. If one has some prior knowledge about some part of the generating system \eqref{eq:generator}, then the results of the paper can be used to easily construct a parameterisation of the estimator \eqref{eq:T0sys}. Secondly, there are many subtleties in the consistency analysis of subspace identification algorithms with inputs \cite{LindquistBook,Katayama:05,chiuso2004ill}, so an alternative approach involving the identification of
an autonomous model \eqref{eq:generator} of $(\y,\w)$ could be advantageous in some cases, i.e., estimating the data generator \eqref{eq:generator} via autonomous system identification, and then using the results of the paper to obtain the estimate of $\y$, is more consistent. That is, with the same amount of data, the estimated predictor will have better performance, in the sense of lower validation mean square error.

\section{Computational Example}
The following examples' code is available on GitLab \cite{gitlab}.
To illustrate the findings consider the system
\begin{subequations}\label{eq:exSystem1}
    \begin{align}       
        \x(t+1)&=\underset{A}{\underbrace{\begin{bmatrix} 1.08 & -0.23 \\ 0.58 & 0.27\end{bmatrix}}}\x(t)+\underset{B}{\underbrace{\begin{bmatrix}-0.56& -1.4 \\-0.56 & -0.6\end{bmatrix}}} v(t),\\
        \begin{bmatrix} \y(t) \\ \w(t) \end{bmatrix} &= \underset{C}{\underbrace{\begin{bmatrix} -0.25& 2.25 \\ 1.24 & -1.25\end{bmatrix}}}\x(t)+\underset{D}{\underbrace{\begin{bmatrix}-0.14& -1 \\ 0 &  -1\end{bmatrix}}} v(t),\\
        v(t)&\sim \mathcal{N}(0,I),
    \end{align}
\end{subequations}
The system \eqref{eq:exSystem1} is such that $\w$ is feedback free from $\y$, later we will find the upper block diagonal form of the innovation process of \eqref{eq:exSystem1}. Before that, we first find the forward innovation form of system \eqref{eq:exSystem1} by following \cite[Chapter 6]{LindquistBook}, for the sake of completeness we reiterate statements of \cite{LindquistBook}. Note \cite{LindquistBook} assumes $v(t)$ in \eqref{eq:exSystem1} to be normalised Gaussian.\\
First find $P$, which solves Lyapunov equation
\begin{align}
    P=APA^T+BB^T, 
\end{align}
this implies
\begin{align}
    P=\begin{bmatrix} 11.34 & 9.22\\ 9.22 & 7.96 \end{bmatrix}
\end{align}
now compute
\begin{align}
    \bar{C}=CPA+DB^T=\begin{bmatrix} 17.24 & 15.28\\ 3.79 & 2.47 \end{bmatrix}\\
    \Lambda_0=CPC^T+DD^T = \begin{bmatrix} 31.64 & 3.72\\ 3.72 & 2.28 \end{bmatrix}
\end{align}
then find $\Pi$, which solves the following algebraic Ricatti equation
\begin{align}
    \Pi&=A\Pi A^T+(\bar{C}^T-A\Pi C^T)\Delta(\Pi )^{-1}(\bar{C}^T-A\Pi C^T)^T,
\end{align}
this implies
\begin{align}
    \Pi&=\begin{bmatrix} 11.1 & 8.98\\ 8.98 & 7.71 \end{bmatrix} \nonumber
\end{align}
with $\Delta(\Pi )=\Lambda_0-C\Pi C^T=\bE[\e_g(t)\e_g^T(t)]=\begin{bmatrix} 2 & 1\\ 1 & 1 \end{bmatrix} $ and finally we compute the gain 
\begin{align}
    K_g=(\bar{C}^T-A\Pi C^T)\Delta(\Pi)^{-1}=\begin{bmatrix} 0.5 & 0.9\\ 0.49 & 0.11 \end{bmatrix},
\end{align}
with which we obtain the system \eqref{eq:exSystem1} in forward innovation form
\begin{subequations} \label{eq:numEx:initInnovSys}
    \begin{align}       
\x(t+1)&=\begin{bmatrix} 1.08 & -0.23\\ 0.58 & 0.27 \end{bmatrix}\x(t)+\begin{bmatrix} 0.5 & 0.9\\ 0.49 & 0.11 \end{bmatrix}\e(t)\\ 
 \begin{bmatrix} \y(t) \\ \w(t) \end{bmatrix} &= \begin{bmatrix} -0.25 & 2.25\\ 1.24 & -1.25 \end{bmatrix} \x(t)+\e(t),\quad Q=\begin{bmatrix} 2 & 1\\ 1 & 1 \end{bmatrix} 
    \end{align}
\end{subequations}
The triangular form is obtained by applying SVD to the observability matrix of $(A_g,C_\w)$. From SVD, contrary to standard practice, we sort the singular values from lowest to highest, and appropriately sort the columns of the matrices $U$ and $V$ containing the left and right singular vectors respectively.  
Using the transformation $$T=V^T=\begin{bmatrix} -0.71 & -0.7\\ -0.7 & 0.71 \end{bmatrix},$$  we obtain a system in triangular form as \eqref{eq:similar}.
\begin{subequations} \label{eq:numEx:diagInnovSys}
    \begin{align}       
        \x(t+1)&=\begin{bmatrix} 0.85 & 0.81\\ 0 & 0.5 \end{bmatrix}\x(t)+\begin{bmatrix} -0.7 & -0.71\\ 0 & -0.56 \end{bmatrix}\e(t)\\ 
 \begin{bmatrix} \y(t) \\ \w(t) \end{bmatrix} &= \begin{bmatrix} -1.41 & 1.77\\ 0 & -1.76 \end{bmatrix} \x(t)+\e(t)
    \end{align}
\end{subequations}
Since a transformation that maps system \eqref{eq:numEx:initInnovSys} to upper block diagonal system \eqref{eq:numEx:diagInnovSys} exists, $\w$ is feedback free from $\y$, and we can compute the estimator of the form \eqref{eq:T0sys}
\begin{subequations} \label{eq:NumEx:pred}
	\begin{align}
\hat{\x}(t+1)&=\begin{bmatrix} 0.85 & -1.69\\ 0 & -0.49 \end{bmatrix}\hat{\x}(t)+\begin{bmatrix} -1.42\\ -0.56 \end{bmatrix}\w(t)\\ 
 \hat{\y}(t) &= \begin{bmatrix} -1.41 & 3.53 \end{bmatrix} \hat{\x}(t)+\w(t)
 \end{align}
\end{subequations}
The system \eqref{eq:NumEx:pred} produces the least square estimate of $\y(t)$, at least when the generating system is fully known. Furthermore, the results of the paper are also useful for system identification. 
\section{System identification example}
The proposed construction could be useful in system identification for parametric methods, since apriori knowledge of the generating system \eqref{eq:generator} can easily be translated to information about the estimator. For example, say we know $(A_g,K_g,C_g,Q)$, except for one element of $A_g$, i.e., $A_{g_{1,2}}=\theta$, then the parameterised generator is given by 
\begin{subequations} \label{eq:numEx:parGen}
    \begin{align}       
\hat{\x}(t+1)&=\begin{bmatrix} 0.85 & \theta \\ 0 & 0.5 \end{bmatrix}\hat{\x}(t)+\begin{bmatrix} -1.41\\ -0.56 \end{bmatrix}\w(t)\\ 
 \hat{\y}(t) &= \begin{bmatrix} -1.41 & 3.52 \end{bmatrix} \hat{\x}(t)+\w(t)
    \end{align}
\end{subequations}
and the estimator is parameterised by
\begin{subequations} 
	\begin{align}
 \hat{\x}(t+1|\theta)&=\begin{bmatrix} 0.85 & \theta -2.49\\ 0 & -0.49 \end{bmatrix}\hat{\x}(t|\theta)+\begin{bmatrix} -1.41\\ -0.56 \end{bmatrix}\w(t)\\ 
 \hat{\y}(t|\theta) &= \begin{bmatrix} -1.41 & 3.52 \end{bmatrix} \hat{\x}(t|\theta)+\w(t)
	\end{align}
\end{subequations}
Then $\theta$ can be found by minimising the mean square error (MSE) on some collected data, i.e., $$\theta^*=\arg\min_\theta (\frac{1}{N}\sum_{t=0}^N(\y(t)(\omega)-\hat{\y}(t\mid \theta))^2).$$ With the results of this paper one can also obtain an estimator of $\y$ from data, by using autonomous system identification on generator system \eqref{eq:generator}, then use the proposed construction to obtain an estimator. Using $100$ points for identification, we obtain validation MSE of $4.88$ when identifying the predictor, which constitutes $78.08\%$ "Variance Accounted For" (VAF)$$\text{VAF}=\min \left \{0, 1-\frac{\text{var}(\y-\hat{\y})}{\text{var}(\y)}\right \}\times 100\%,$$ and when identifying the generator we obtain validation MSE of $4.69$ ($78.61\%$ VAF). \par 
To better explore the better consistency of first identifying the 'data generator', i.e. the system in forward innovation form \ref{as:generator}, and then using the proposed construction, let us consider a larger system. That is, a system with $n=10$ states, of which $p_2=6$ correspond to process $\w$, i.e. $A_{2,2}\in \mathbb{R}^{6\times6}$ and $p_1=4$, to process $\y$ i.e. $A_{1,1}\in \mathbb{R}^{4\times 4}$, with $p=3$ outputs, i.e. $\y(t)\in\mathbb{R}^3$, and $q=2$, i.e. $\w(t)\in \mathbb{R}^2$. We will consider 4 cases for system identification of this $n=10$ state system.
\begin{enumerate}[label=Case \arabic*), leftmargin=*]
    \item We will assume no prior information, and the identification will consist of identifying the full system, by minimising MSE.
    \begin{enumerate}[label=Case 1.\arabic*]
        \item \textbf{"No prior information, estimate predictor"} \\
        Identification of the system in the form of the \emph{estimator} \eqref{eq:T0sys}
        \item \textbf{"No prior information, est. generator"} \\
        Identification of the system in the form of the \emph{generator} \eqref{eq:generator}, and after identification compute the 'optimal' estimator as described in Section \ref{sec:Result}
    \end{enumerate}
    \item We will assume some prior information, more specifically we will assume we are given $(A_{2,2},K_{2,2},C_{2,2},Q_{2,2})$, and the rest of the system is unknown, and thus parameterised, i.e.,
    Note that all parameterised matrices $A_{1,1}(\theta),A_{1,2}(\theta),K_{1,1}(\theta),K_{1,2}(\theta),C_{1,1}(\theta), C_{1,2}(\theta),Q_{1,2}(\theta)$ are fully parameterised, that is each element of the matrix is assigned a parameter.
    \begin{enumerate}[label=Case 2.\arabic*]
        \item \textbf{"Some prior information, estimate predictor"}\\
        Identification of the system in the form of the \emph{estimator} \eqref{eq:T0sys}
        \item \textbf{"Some prior information, estimate generator"}\\
        Identification of the system in the form of the \emph{generator} \eqref{eq:generator}, and after identification compute the 'optimal' estimator as described in Section \ref{sec:Result}
    \end{enumerate}
    \item[Case 0] \textbf{"Full knowledge of the system"}\\
    For comparison, we will include statistics in Table \ref{tab:stats}, of the optimal estimator, obtained from the results of the paper.
\end{enumerate}
    \begin{subequations}
            \begin{align}
        \tilde{\x}(t+1)&=\begin{bmatrix} A_{1,1}(\theta) & A_{1,2}(\theta) \\ 0 & A_{2,2} \end{bmatrix}\tilde{x}(t)+\begin{bmatrix} K_{1,1}(\theta) & K_{1,2}(\theta) \\ 0 & K_{2,2}\end{bmatrix}\e(t)\\
        \begin{bmatrix} \tilde{\y} \\ \tilde{\w} \end{bmatrix} &= \begin{bmatrix} C_{1,1}(\theta) & C_{1,2}(\theta)\\0 & C_{2,2} \end{bmatrix} \tilde{\x}(t) + \e(t)\\
        \e(t)&\sim\mathcal{N}\left ( 0 , \begin{bmatrix} Q_{1,1}& Q_{1,2}(\theta)\\ Q_{1,2}^T(\theta) & Q_{2,2} \end{bmatrix} \right )
    \end{align}
    \end{subequations}
All system identification is done in \textit{Matlab}, using "idss" models for Case 1 and "idgrey" models for Case 2, and prediction error is minimised using \textit{Matlab}'s "pem" function.

In order to compare the different cases and approaches to identifying the 'best' estimator of $\y$ from data, we will consider 2 indicators: average validation MSE, and average VAF. The average is understood 
If we collect $M$ trajectories $\{\{\y(t)(\omega_m),\w(t)(\omega_m)\}_{t\in\{1,\dots,N\}}\}_{m=\{1,\dots,M\}}$, then for each trajectory we perform system identification, which minimises $(\text{MSE})_m$, i.e.
\begin{align*}
    \text{(MSE)}_m&=\frac{1}{N}\sum_{t=1}^N \tilde{y}^T(t,m)\tilde{y}(t,m),\\
    \tilde{y}(t,m)&=\y(t)(\omega_m)-\hat{\y}(t)(\omega_m).
\end{align*}
Then the average MSE, is given by
\begin{align}
    \text{average MSE}=\frac{1}{M}\sum_{m=1}^M \text{(MSE)}_m \label{eq:averageStat}
\end{align}

Average VAF should be understood in a similar manner. The difference between average MSE and average validation MSE, is that we obtain the estimator by minimising MSE on one set of data, and achieved MSE is labeled training MSE. Whereas, for validation MSE, we take a previously identified estimator and compute MSE on a different data set. Note that in Table \ref{tab:stats} tra\par 
\begin{figure}[h]
\centering
\includegraphics[width=\linewidth]{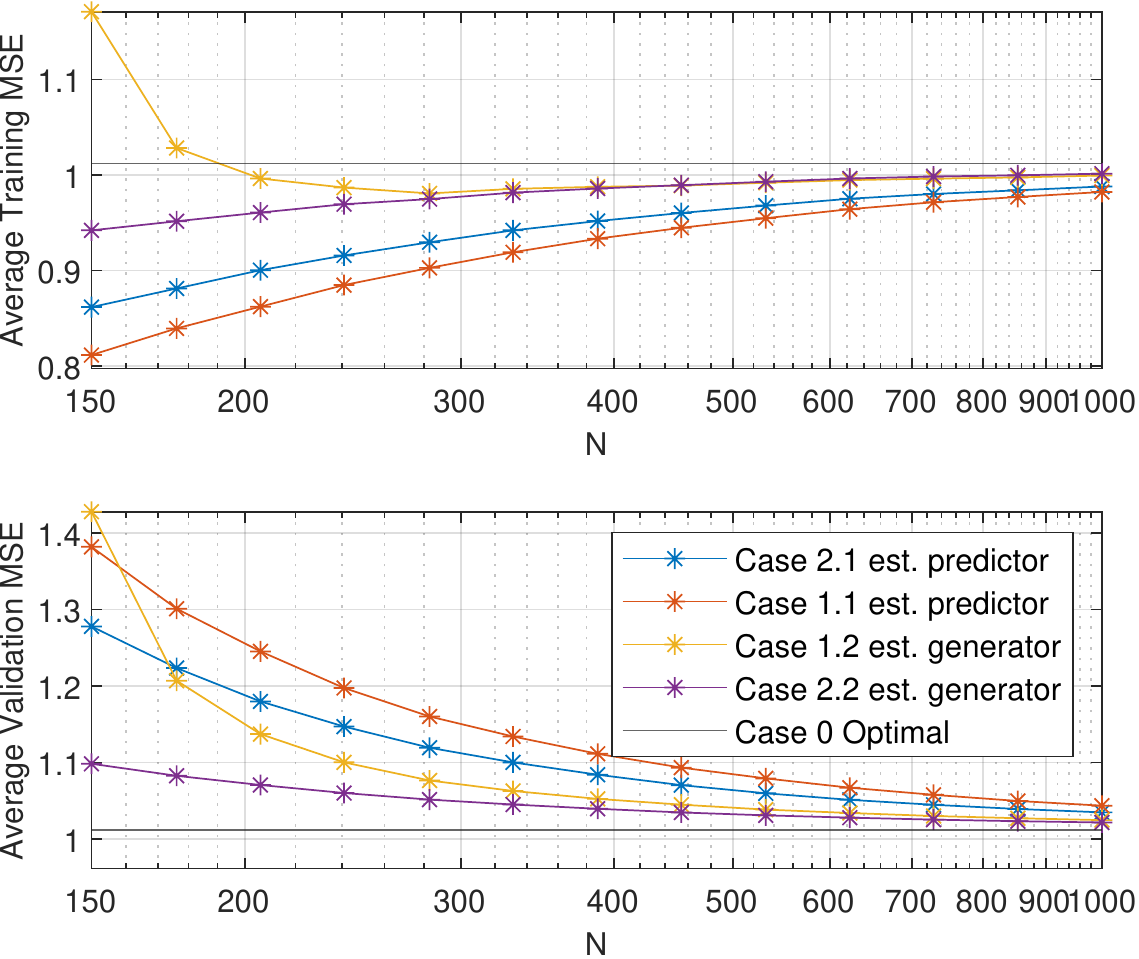} 
\caption{For a given $N$, Multiple trajectories have been generated, so we have a set of $M$ data sets $S_{N,i}$ with $N$ samples each, then we try to estimate the predictor from each data set $S_{N,i}$ using the 4 approaches discussed above, and compute the average training MSE and average validation MSE. Note: training MSE is computed only for those $N$ data points that have been used to estimate the predictor, whereas validation MSE shows the true predictive power of the predictor. } 
\label{fig:Average}
\end{figure}
\begin{table*}[ht]
    \centering
\caption{Summary of the results. Best values marked in bold. All quantities are averaged over multiple system identifications with varying data, as explained in \eqref{eq:averageStat}}
\begin{tabular}{rr|p{1.5cm}|p{2.1cm}|p{2.1cm}|p{2.1cm}|p{2.1cm}|}
& & Case 0 Optimal & Case 1.1 est. predictor & Case 1.2 est. generator & Case 2.1 est. predictor & Case 2.2 est. generator \\ 
\hline 
\multicolumn{2}{r|}{Total Parameters} & 0 & 156 & 150 & 96 & 90 \\ 
\hline \multicolumn{1}{r|}{\multirow{6}{*}{\centering $N=150$}} & Training MSE & — & \textbf{0.812} & 1.171 & 0.862 & 0.942 \\ 
\multicolumn{1}{r|}{} & Validation MSE & 1.012 & 1.382 & 1.428 & 1.278 & \textbf{1.098} \\ 
\multicolumn{1}{r|}{} & Validation $VAF_1$ & 22.7\% & 5.6\% & 11.6\% & 8.5\% & \textbf{16.6}\% \\ 
\multicolumn{1}{r|}{} & Validation $VAF_2$ & 81\% & 73.2\% & 76\% & 74.9\% & \textbf{79.6}\% \\ 
\multicolumn{1}{r|}{} & Validation $VAF_3$ & 92.6\% & 90.6\% & 89.5\% & 91\% & \textbf{91.9}\% \\ 
\multicolumn{1}{r|}{} & Validation $\overline{VAF}$ & 65.4\% & 56.5\% & 59.1\% & 58.2\% & \textbf{62.7}\% \\ 
\hline \multicolumn{1}{r|}{\multirow{6}{*}{\centering $N=1000$}} & Training MSE & — & \textbf{0.982} & 1 & 0.988 & 1.002 \\ 
\multicolumn{1}{r|}{} & Validation MSE & 1.012 & 1.044 & 1.025 & 1.035 & \textbf{1.022} \\ 
\multicolumn{1}{r|}{} & Validation $VAF_1$ & 23\% & 20.6\% & 22\% & 20.9\% & \textbf{22.3}\% \\ 
\multicolumn{1}{r|}{} & Validation $VAF_2$ & 81.7\% & 81\% & \textbf{81.5}\% & 81.1\% & \textbf{81.5}\% \\ 
\multicolumn{1}{r|}{} & Validation $VAF_3$ & 92.7\% & 92.4\% & 92.5\% & 92.5\% & \textbf{92.6}\% \\ 
\multicolumn{1}{r|}{} & Validation $\overline{VAF}$ & 65.8\% & 64.7\% & 65.3\% & 64.8\% & \textbf{65.5}\% \\ 
\hline 
\end{tabular}
    \label{tab:stats}
\end{table*}
In Figure \ref{fig:Average}, we see how on average the number of samples used for training affect the training MSE and validation MSE. Something of note is the high validation MSE of Case 1.2 identifying the full generator, note that this case has the highest number of parameters (156) to estimate as indicated in Table \ref{tab:stats}, and thus is \emph{over-parameterised}. If we ignore the over-parameterisation, we see that on average we are better off identifying the generator, instead of identifying the estimator.  \par

Since the system in this example $\y$ has 3 components, Table \ref{tab:stats} reports "Variance Accounted For" for each component, i.e. with $\y_{i,m},\hat{\y}_{i,m}$ as short-hand for the $i$th component of trajectories $\y(\omega_m)$ and appropriately $\hat{\y}(\omega_m)$.
\begin{align*}
    \text{(VAF)}_{m,i}&=\max \left \{0, 1-\frac{var(\y_{i,m}-\hat{\y}_{i,m})}{var(\y_{i,m})} \right \} \times 100\%\\
    \text{(VAF)}_{i}&=\frac{1}{M}\sum_{m=1}^M\text{(VAF)}_{m,i}
\end{align*}
this shows how well each component of $\y$ is estimated, for example using $150$ samples for identification we can at best estimate $16.6\%$ of the first component's variance. Whereas, if we had known the system (Case 0) we should be able to estimate $22.7\%$ of the first component. If, in this hypothetical example, we had collected $N=1000$ samples for identification, then we could estimate $22.3\%$ of the first component, much closer to best possible estimation (Case 0).\\
We also report the average over components VAF, i.e.
\begin{align*}
    \overline{\text{VAF}}=\frac{1}{p}\sum_{i=1}^p \text{(VAF)}_i
\end{align*}
In summary if we collect enough samples to avoid over-parameterisation, then identifying the generator, and then applying the construction as defined in Section \ref{sec:Result}, yields better performing estimators of $\y$.




\bibliographystyle{plain}        
\bibliography{bib}           



\end{document}